\documentclass{amsart}
\usepackage{amssymb}
\usepackage{amsfonts}
\usepackage{amsmath}
\usepackage{graphicx}
\usepackage[pagebackref]{hyperref}
\theoremstyle{plain}
\newtheorem{thm}{Theorem}[section]
\newtheorem{corollary}[thm]{Corollary}
\newtheorem{lemma}[thm]{Lemma}

\theoremstyle{definition}

\theoremstyle{remark}

\newtheorem{example}[thm]{Example}

\newtheorem{claim}{\rm Claim}

 \DeclareMathOperator{\qf}{qf}
 
 \DeclareMathOperator{\Ker}{Ker}
\DeclareMathOperator{\Ima}{Im} 

\DeclareMathOperator{\Spec}{Spec} \DeclareMathOperator{\Max}{Max}
\DeclareMathOperator{\Cl}{Cl} \DeclareMathOperator{\Pic}{Pic}
\numberwithin{equation}{section}

\def\1{{\rm (1)}}
\def\2{{\rm (2)}}
\def\3{{\rm (3)}}
\def\4{{\rm (4)}}
\def\5{{\rm (5)}}

\begin{document}

\title[Constituent groups of Clifford semigroups]{Constituent groups of Clifford semigroups \\arising from $t$-closure}

\author{S. Kabbaj}

\address{Department of Mathematical Sciences, King Fahd University of Petroleum \& Minerals, P.O. Box 5046, Dhahran 31261, KSA}

\email{kabbaj@kfupm.edu.sa}

\author{A. Mimouni}

\address{Department of Mathematical Sciences, King Fahd University of Petroleum \& Minerals, P.O. Box 5046, Dhahran 31261, KSA}

\email{amimouni@kfupm.edu.sa}

\thanks{This work was funded by KFUPM under Project \# MS/t-Class/257.}

\dedicatory{\small\tt to Alain Bouvier on the occasion of his sixty-fifth 
birthday}

\subjclass[2000]{13C20, 13F05, 11R65, 11R29, 20M14}

\keywords{Class semigroup, $t$-class
semigroup, class group, $t$-ideal, $t$-closure, Clifford semigroup, Picard group, valuation
domain, Pr\"ufer domain, Krull-type domain, PVMD}

\begin{abstract}
The $t$-class semigroup of an integral domain $R$, denoted ${\mathcal S}_{t}(R)$, is the semigroup
of fractional $t$-ideals modulo its subsemigroup of nonzero principal ideals with the operation
induced by ideal $t$-multiplication. We recently proved that if $R$ is a Krull-type domain (in the
sense of Griffin), then ${\mathcal S_t}(R)$ is a Clifford semigroup. This paper aims to describe
the idempotents of ${\mathcal S_t}(R)$ and the structure of their associated groups. We extend and
recover well-known results on class semigroups of valuation domains and Pr\"ufer domains of finite
character.
\end{abstract}

\maketitle

\section{Introduction}
The $t$-operation in integral domains is considered as one of the keystones of multiplicative ideal
theory. It originated in Jaffard's 1960 book ``Les Syst\`emes d'Id\'eaux" \cite{J} and was
investigated by many authors in the 1980's. From the $t$-operation stemmed the notion of class
group of an arbitrary domain, extending both notions of divisor class group (in Krull domains) and
Picard group (in Pr\"ufer domains). Class groups were first introduced and developed by Bouvier
\cite{B} and Bouvier \& Zafrullah \cite{BZ}, and have been, since then, extensively studied  in the
literature. In the 1990's, the attention of some authors moved from the class group to the class
semigroup, first considering orders in number fields \cite{ZZ} and valuation domains \cite{BS}, and
then more general contexts \cite{Ba1,Ba2,Ba3,Ba4,KM1}. The basic idea is to look at those domains
that have Clifford class semigroup.

Let $R$ be an integral domain with quotient field $K$. For a nonzero fractional ideal $I$ of $R$,
let $I^{-1}:=(R:I)=\{x\in K\mid xI\subseteq R\}$. The $v$- and $t$-closures of $I$ are defined,
respectively, by $I_v:=(I^{-1})^{-1}$ and $I_t:=\cup J_v$ where $J$ ranges over the set of finitely
generated subideals of $I$. The ideal $I$ is said to be a $v$-ideal if $I_v=I$, and a $t$-ideal if
$I_t=I$. Under the ideal $t$-multiplication $(I,J)\mapsto (IJ)_t$, the set $F_{t}(R)$ of fractional
$t$~-ideals of $R$ is a semigroup with unit $R$. An invertible element for this operation is called
a $t$-invertible $t$-ideal of $R$. So the set $Inv_{t}(R)$ of $t$-invertible fractional
$t$-ideals of $R$ is a group with unit $R$ (Cf. \cite{Gi}). Let $F(R)$, $Inv(R)$, and $P(R)$ denote
the sets of nonzero, invertible, and nonzero principal fractional ideals of $R$, respectively.
Under this notation, the Picard group \cite{A,BM,G}, class group \cite{B,BZ}, $t$-class semigroup
\cite{KM2}, and class semigroup \cite{BS,KM1,ZZ} of $R$ are defined as follows:
$$\Pic(R):=\frac{Inv(R)}{P(R)}\subseteq \Cl(R):=\frac{Inv_{t}(R)}{P(R)}\subseteq
S_{t}(R):=\frac{F_{t}(R)}{P(R)}\subseteq S(R):=\frac{F(R)}{P(R)}.$$

A commutative semigroup $S$ is said to be Clifford if every element
$x$ of $S$ is (von Neumann) regular, i.e., there exists $a\in S$
such that $x^{2}a=x$. The importance of a Clifford semigroup $S$
resides in its ability to stand as a disjoint union of subgroups
$G_e$, where $e$ ranges over the set of idempotents of $S$, and
$G_e$ is the largest subgroup of $S$ with identity equal to $e$; namely, $G_e=\{ae\mid abe=e$ for some $b\in S\}$ (Cf.
\cite{Ho}). Often, the $G_e$'s are called the constituent groups of $S$.

A domain $R$ is called a PVMD (Pr\"ufe $v$-multiplication domain) if
$R_{M}$ is a valuation domain for each $t$-maximal ideal $M$ of $R$.
Ideal $t$-multiplication converts ring notions such as PID,
Dedekind, Bezout, and Pr\"ufer, respectively to UFD, Krull, GCD, and
PVMD. Recall at this point that the PVMDs of finite $t$-character
(i.e., each proper $t$-ideal is contained in only finitely many
$t$-maximal ideals) are exactly the Krull-type domains introduced by
Griffin in 1967-68 \cite{Gr1,Gr2}.

Divisibility properties of $R$ are often reflected in semigroup-theoretic properties of $S(R)$ and
$S_{t}(R)$. Obviously, Dedekind (resp., Krull) domains have Clifford class (resp., $t$-class)
semigroup. In 1994, Zanardo and Zannier proved that all orders in quadratic fields have Clifford
class semigroup \cite{ZZ}. They also showed that the ring of all entire functions in the complex
plane (which is Bezout) fails to have this property. In 1996, Bazzoni and Salce proved that any
arbitrary valuation domain has Clifford semigroup \cite{BS}. In \cite{Ba1,Ba2,Ba3}, Bazzoni
examined the case of Pr\"ufer domains of finite character, showing that these, too, have Clifford
class semigroup. In 2001, she completely resolved the problem for the class of integrally closed
domains by stating that ``{\em an integrally closed domain has Clifford class semigroup if and only
if it is a Pr\"ufer domain of finite character}" \cite[Theorem 4.5]{Ba4}. Recently, we extended
this result to the class of PVMDs of finite character; namely, ``{\em a PVMD has Clifford $t$-class
semigroup if and only if it is a Krull-type domain}" \cite[Theorem 3.2]{KM2}.

This paper extends Bazzoni and Salce's study of groups in the class semigroup of a valuation domain
\cite{BS} or a Pr\"ufer domain of finite character \cite{Ba2,Ba3} to a larger class of integral
domains. Precisely, we describe the idempotents of $S_{t}(R)$ and the structure of their associated
groups when $R$ is a Krull-type domain. Indeed, we prove that there are two types of idempotents in
${\mathcal S_t}(R)$: those represented by certain fractional overrings of $R$ and those represented
by finite intersections of $t$-maximal ideals of certain fractional overrings of $R$. Further, we
show that the group associated with an idempotent of the first type equals the class group of the
fractional overring, and characterize the elements of the group associated with an idempotent of
the second type in terms of their localizations at $t$-prime ideals. Our findings recover Bazzoni's
results on the constituents groups of the class semigroup of a Pr\"ufer domain of finite character.

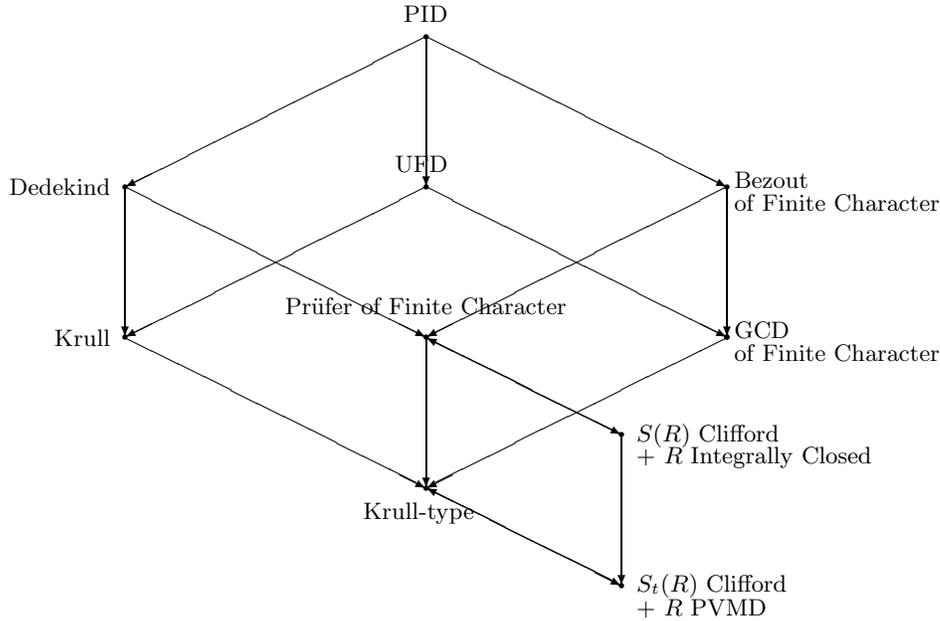
\begin{figure}[h!]
\centering
\[\setlength{\unitlength}{1mm}
\begin{picture}(100,72)(0,-100)
\put(40,-30){\vector(0,-1){20}}\put(40,-30){\vector(2,-1){40}}
\put(40,-30){\vector(-2,-1){40}}\put(0,-50){\vector(0,-1){20}}
\put(0,-50){\vector(2,-1){40}}
\put(80,-50){\vector(0,-1){20}}\put(80,-50){\vector(-2,-1){40}}
 \put(66,-83){\vector(0,-1){20}}
 \put(40,-50){\vector(2,-1){40}}
\put(40,-50){\vector(-2,-1){40}}\put(0,-70){\vector(2,-1){40}}
\put(80,-70){\vector(-2,-1){40}} \put(40,-70){\vector(0,-1){20}}
\put(40,-70){\vector(2,-1){26}}\put(40,-90){\vector(2,-1){26}}
\put(66,-103){\vector(-2,1){26}}\put(66,-83){\vector(-2,1){26}}
\put(40,-30){\circle*{.7}} \put(40,-28){\makebox(0,0)[b]{\small
PID}} \put(0,-50){\circle*{.7}} \put(-2,-50){\makebox(0,0)[r]{\small
Dedekind}} \put(40,-50){\circle*{.7}}
\put(39.4,-48){\makebox(0,0)[b]{\small UFD}}
\put(80,-50){\circle*{.7}} \put(81,-49){\makebox(0,0)[l]{\small
Bezout }}\put(81,-52){\makebox(0,0)[l]{\small of Finite Character}}
\put(0,-70){\circle*{.7}} \put(-2,-70){\makebox(0,0)[r]{\small
Krull}} \put(80,-70){\circle*{.7}}
\put(81,-69){\makebox(0,0)[l]{\small GCD }}
\put(80,-70){\circle*{.7}} \put(81,-72){\makebox(0,0)[l]{\small of
Finite Character}}
\put(40,-70){\circle*{.7}}\put(40,-67){\makebox(0,0)[b]{\small
Pr\"ufer of Finite Character}}
 \put(40,-90){\circle*{.7}} \put(39,-92){\makebox(0,0)[t]{\small Krull-type}}
\put(66,-83){\circle*{.7}}\put(68,-83){\makebox(0,0)[l]{\small $S(R)$
Clifford}}\put(68,-86){\makebox(0,0)[l]{\small + $R$ Integrally Closed}}
\put(66,-103){\circle*{.7}}\put(68,-103){\makebox(0,0)[l]{\small $S_{t}(R)$
Clifford}}\put(68,-106){\makebox(0,0)[l]{\small + $R$ PVMD}}
\end{picture}\]\medskip

\caption{\tt A ring-theoretic perspective for Clifford
property}\label{D1}
\end{figure}

All rings considered in this paper are integral domains. For the
convenience of the reader, Figure~\ref{D1} displays a
 diagram of implications summarizing the relations between the main
classes of integral domains involved in this work. It also places
the Clifford property in a ring-theoretic perspective.

\section{Main Result}\label{sec:1}

An overring $T$ of a domain $R$ is $t$-linked over $R$ if
$I^{-1}=R\Rightarrow(T:IT)=T$, for each finitely generated ideal $I$
of $R$ \cite{AHZ,KP}. In Pr\"ufer domains, the $t$-linked property
collapses merely to the notion of overring (since every finitely
generated proper ideal is invertible and then different from $R$).
This concept plays however a crucial role in any attempt to extend
classical results on Pr\"ufer domains to PVMDs (via $t$-closure).
Recall also that an overring $T$ of $R$ is fractional if $T$ is a
fractional ideal of $R$; in this case, any (fractional) ideal of $T$
is a fractional ideal of $R$. Of significant importance too for the
study of $t$-class semigroups is the notion of $t$-idempotence;
namely, a $t$-ideal $I$ of a domain $R$ is $t$-idempotent if
$(I^{2})_{t} = I$.

The following discussion, connected with the $t$-ideal structure of a PVMD, will be of use in the
sequel without explicit mention. Let $R$ be a PVMD. Note that prime ideals of $R$ contained in
a  $t$-maximal ideal are necessarily $t$-ideals and form a chain \cite{Kg}. Also, recall that
$t$-linked overrings of $R$ are exactly the subintersections of $R$; precisely, $T$ is a $t$-linked
overring of $R$ if and only if $T=\bigcap R_{P}$, where $P$ ranges over some set of $t$-prime
ideals of $R$ \cite{CZ,Kg2}. Further, for any $t$-ideal $I$ of $R$, $(I:I)$ is a subintersection
of $R$ and hence a fractional $t$-linked overring of $R$. Indeed, let $\Max_{t}(R, I):=\{M\in
\Max_{t}(R)\ | \ I\subseteq M\}$ and $\overline{\Max_{t}}(R, I):=\{M\in \Max_{t}(R)\ |\ I\nsubseteq
M\}$, where $\Max_{t}(R)$ denotes the set of
$t$-maximal ideals of $R$. Then $(I:I)=(\bigcap_{\alpha}R_{M_{\alpha}})\cap(\bigcap_{\beta}R_{N_{\beta}})$, where $M_\alpha$
ranges over $\overline{\Max_{t}}(R, I)$ and $N_{\beta}$ denotes the set of zerodivisors of
$R_{M_{\beta}}$ modulo $IR_{M_{\beta}}$ where $M_{\beta}$ ranges over $\Max_{t}(R,I)$. Consequently,
$(I:I)$ is a PVMD \cite{Kg}. Finally, given $M_{1}$ and $M_{2}$ two $t$-maximal ideals of $R$,
we will denote by $M_{1}\wedge M_{2}$ the largest prime ideal of $R$ contained in $M_{1}\cap M_{2}$.

Throughout, we shall use  ${\overline I}$ to denote the isomorphy class of an ideal $I$ of $R$ in
$S_{t}(R)$, $\qf(R)$ to denote the quotient field of $R$. Recall that the class group of an integral domain $R$, denoted $\Cl(R)$,
is the group of fractional $t$-invertible $t$-ideals modulo its subgroup of nonzero principal
fractional ideals. Also we shall use $v_{1}$ and $t_{1}$ to denote the $v$- and $t$-operations with
respect to an overring $T$ of $R$. By \cite[Theorem 3.2]{KM2}, if $R$ is a Krull-type domain, then $S_{t}(R)$ is Clifford
and hence  a disjoint union of subgroups $G_{\overline{J}}$, where $\overline{J}$ ranges over the
set of idempotents of $S_{t}(R)$ and $G_{\overline{J}}$ is the largest subgroup of $S_{t}(R)$ with
unit $\overline{J}$. Notice for convenience that in valuation and Pr\"ufer domains the $t$- and
trivial operations (and hence the $t$-class and class semigroups) coincide. At this point, it is
worthwhile recalling Bazzoni-Salce result that valuation domains have Clifford class semigroup
\cite{BS}. Next we announce the main result of this paper.

\begin{thm}\label{sec:1.2} Let $R$ be a Krull-type domain and
$I$ a $t$-ideal of $R$. Set $T:=(I:I)$ and $\Gamma(I) :=\{$finite intersections of
$t$-idempotent $t$-maximal ideals of $T\}$. Then:\\ $\overline{I}$ is an idempotent of $S_{t}(R)$
if and only if there exists a unique $J\in \{T\}\cup \Gamma(I)$ such
that $\overline{I}=\overline{J}$. Moreover, \\
\1 if $J=T$, then $G_{\overline{J}}\cong\Cl(T)$;\\
\2 if $J=\bigcap_{1\leq i\leq r}Q_{i}\in\Gamma(I)$, then the
following sequence of natural group homomorphisms is exact
$$0\longrightarrow \Cl(T)\stackrel{\phi}\longrightarrow
G_{\overline{J}}\stackrel{\psi}\longrightarrow \prod_{1\leq i\leq
r}G_{\overline{Q_{i}T_{Q_{i}}}}\longrightarrow 0$$  where $G_{\overline{Q_{i}T_{Q_{i}}}}$ denotes
the constituent group of the Clifford semigroup $S(T_{Q_{i}})$ associated with
$\overline{Q_{i}T_{Q_{i}}}$.
\end{thm}

The proof of the theorem involves several preliminary lemmas, some
of which are of independent interest. We often will be appealing to
some of them without explicit mention.

\begin{lemma}\label{sec:1.0-0}\cite[Lemma 2.1]{KM2}
Let $I$ be a $t$-ideal of a domain $R$. Then
 ${\overline I}$ is regular in $S_{t}(R)$ if and only if
$I=(I^{2}(I:I^{2}))_{t}$. \qed
\end{lemma}

\begin{lemma}\label{sec:1.10}
Let $R$ be a PVMD, $T$ a $t$-linked overring of $R$, and $Q$ a $t$-prime ideal of $T$. Then $P=Q\cap R$ is a $t$-prime ideal of $R$ with $R_{P}=T_{Q}$. If, in addition, $Q$ is supposed to be $t$-idempotent in $T$, then so is $P$ in $R$.
\end{lemma}

\proof
 Since $R$ is a PVMD, by \cite[Proposition 2.10]{KP}, $T$ is $t$-flat over $R$. Hence
$R_{P}=T_{Q}$. Moreover, since $T$ is $t$-linked over $R$, then
$P_{t}\subsetneq R$ \cite[Proposition 2.1]{DHLZ}. Hence $P$ is a
$t$-prime ideal of $R$ \cite[Corollary 2.47]{Kg}. Next assume that $(Q^2)_{t_{1}}=Q$. Then
$P^{2}R_{P}=Q^{2}T_{Q}=(Q^{2})_{t_{1}}T_{Q}=QT_{Q}=PR_{P}$ by
\cite[Lemma 3.3]{KM2}. Now, let $M$ be an arbitrary $t$-maximal ideal
of $R$. We claim that $P^{2}R_{M}=PR_{M}$. Indeed, without loss of
generality we may assume $P\subsetneq M$. So $R_{M}\subseteq R_{P}$
and hence $PR_{M}\subseteq PR_{P}$. If $PR_{M}\subsetneq PR_{P}$ and
$x\in PR_{P}\setminus PR_{M}$, necessarily $PR_{M}\subset xR_{M}$
since $R_{M}$ is a valuation domain. Hence, by \cite[Theorem 3.8 and
Corollary 3.6]{HP}, $x^{-1}\in
(R_{M}:PR_{M})=(PR_{M}:PR_{M})=(R_{M})_{PR_{M}}=R_{P}$, absurd.
Therefore $PR_{M}= PR_{P}$. It follows that $P^{2}R_{M}=
P^{2}R_{P}=PR_{P}=PR_{M}$. By \cite[Theorem 2.9]{Kg}, $(P^2)_{t}=P$,
as desired.\qed

\begin{lemma}\label{sec:1.0-1} Let $R$ be a PVMD and $T$ a
$t$-linked overring of $R$. Let $J$ be a common (fractional) ideal of $R$ and $T$. Then the following assertions hold:\\
\1 $J_{t_{1}}=J_{t}$.\\
\2 $J$ is a $t$-idempotent $t$-ideal of $R$ if and only if  $J$ is a
$t$-idempotent $t$-ideal of $T$.
\end{lemma}

\proof \1 Let $x\in J_{t_{1}}$. Then there exists a finitely generated
ideal $B:=\sum_{1\leq i\leq n}a_{i}T$ of $T$ such that $B\subseteq
J$ and $x(T:B)\subseteq T$. Clearly, $A:=\sum_{1\leq i\leq n}a_{i}R$
is a finitely generated ideal of $R$ with $AT=B$. Therefore
$(R:A)\subseteq (T:B)$ and hence $xA(R:A)\subseteq xB(T:B)\subseteq
B\subseteq J$. Moreover $A$ is $t$-invertible in $R$ since $R$ is a
PVMD. It follows that $xR=x(A(R:A))_{t}=(xA(R:A))_{t}\subseteq
J_{t}$. Hence $J_{t_{1}}\subseteq J_{t}$. Conversely, let $x\in
J_{t}$. Then there exists a finitely generated subideal $A$ of $J$
such that $x(R:A)\subseteq R$. Let $N\in\Max_{t}(T)$ with $M:=N\cap
R$. So $x(AT_{N})^{-1}=x(AR_{M})^{-1}=xA^{-1}R_{M}\subseteq
R_{M}=T_{N}$. Therefore $x$ lies in the $v$-closure of $AT_{N}$ in
the valuation domain $T_{N}$. Further $AT_{N}$ is principal (since
it is finitely generated), hence a $v$-ideal. So that $x\in
AT_{N}\subseteq JT_{N}$. Consequently, $x\in
\bigcap_{N\in\Max_{t}(T)}JT_{N}=J_{t_{1}}$, as desired.\\
\2 Straightforward via (1).

\begin{lemma}\label{sec:1.0-2} Let $R$ be a PVMD, $I$ a $t$-ideal of
$R$, and $T:=(I:I)$. Let $J:=\bigcap _{1\leq i\leq r}Q_{i}$, where
each $Q_{i}$ is a $t$-maximal ideal of $T$. Then $J$ is a fractional
$t$-idempotent $t$-ideal of $R$.
\end{lemma}

\proof
Clearly, $J$ is a $t$-ideal of $T$ and hence a fractional $t$-ideal of $R$ by
Lemma~\ref{sec:1.0-1}, with $J=(\prod_{1\leq i\leq r}Q_{i})_{t_{1}}= (\prod_{1\leq i\leq
r}Q_{i})_{t}$. Further, $J^{2}$ is a common ideal of $R$ and $T$, so that
$(J^{2})_{t}=(J^{2})_{t_{1}}=(\prod_{1\leq i\leq r}Q_{i}^{2})_{t_{1}}= (\prod_{1\leq i\leq
r}(Q_{i}^{2})_{t_{1}})_{t_{1}}=(\prod_{1\leq i\leq r}Q_{i})_{t_{1}}=J$, as desired. \qed

\begin{lemma}\label{sec:1.1} Let $R$ be a PVMD, $I$ a
$t$-idempotent $t$-ideal of $R$, and $M$ a $t$-maximal ideal of $R$
containing $I$. Then $IR_{M}$ is an idempotent (prime) ideal of
$R_{M}$.\end{lemma}

\proof     Let $x\in IR_{M}$. Then $x\mu\in I=(I^{2})_{t}$ for some $\mu\in R\setminus
M$. So there exists a finitely generated ideal $B$ of $R$ such that $B\subseteq I^{2}$ and
$x\mu(R:B)\subseteq R$. Therefore $x(R_{M}:BR_{M})=x(R:B)R_{M}=x\mu(R:B)R_{M}\subseteq R_{M}$.
Hence $x$ lies in the $v$-closure of $BR_{M}$. Similar arguments as in the proof of Lemma~\ref{sec:1.0-1}(2) yields that $BR_{M}$ is a
$v$-ideal and hence $x\in BR_{M}\subseteq I^{2}R_{M}$. Therefore $IR_{M}$ is an idempotent ideal
which is necessarily prime since $R_{M}$ is a valuation domain.\qed

\begin{lemma}\label{sec:1.3-0}
Let $R$ be a Krull-type domain, $L$ a $t$-ideal of $R$, and ${\overline J}$ an idempotent of
$S_{t}(R)$. Then $\overline{L}\in G_{{\overline J}}$ if and only if $(L:L)=(J:J)$ and
$(JL(L:L^{2}))_{t}=(L(L:L^{2}))_{t}=(L(J:L))_{t}=J$.
\end{lemma}

\proof
Suppose $\overline{L}\in G_{{\overline J}}$. Then $J=(LK)_{t}$ for some fractional ideal $K$ of $R$
and $xL=(LJ)_{t}$ for some $0\not=x\in\qf(R)$. Moreover,
$((LJ)_{t}J)_{t}=(L(J^{2})_{t})_{t}=(LJ)_{t}=xL$ and also $((LJ)_{t}J)_{t}=x(LJ)_{t}$. So
$L=(LJ)_{t}$. We have $(L:L)\subseteq (LK:LK)\subseteq ((LK)_{t}:(LK)_{t})=(J:J)$ and
$(J:J)\subseteq (LJ:LJ)\subseteq ((LJ)_{t}:(LJ)_{t})=(L:L)$, so that $(L:L)=(J:J)$. Further,
$(JL(L:L^{2}))_{t}=(L(L:L^{2}))_{t}=(L((J:J):L))_{t}=(L(J:JL))_{t}=(L(J:(JL)_{t}))_{t}=(L(J:L))_{t}$.
Clearly, $(L(J:L))_{t}\subseteq J$. Also $LK\subseteq (LK)_{t} = J$. Then $K\subseteq (J:L)$, hence
$J=(LK)_{t}\subseteq (L(J:L))_{t}$. Conversely, take $K:=J(L:L^{2})$ and notice that
$J=(JL(L:L^{2}))_{t}\subseteq (JL)_{t}\subseteq J$, that is, $J=(JL)_{t}$, completing the proof of
the lemma. \qed

\begin{lemma}\label{sec:1.7}
Let $R$ be a PVMD and $I$ a $t$-ideal of $R$. Then\\
\1 $I$ is a $t$-ideal of $(I:I)$.\\
\2 If $R$ is Clifford $t$-regular, then so is $(I:I)$.
\end{lemma}

\proof \1 $(I:I)$ is a $t$-linked overring of $R$ and then apply Lemma~\ref{sec:1.0-1}(1).\\
\2 Let $J$ be a $t$-ideal of $T:=(I:I)$. By Lemma~\ref{sec:1.0-1}(1), $J$ is a $t$-ideal of
$R$. Next, assume that $R$ is Clifford $t$-regular. By \cite[Lemma 3.3]{KM2},
$(J^{2}(J:J^{2}))_{t_{1}}T_{N}=(J^{2}(J:J^{2}))T_{N}=(J^{2}(J:J^{2}))R_{M}=(J^{2}(J:J^{2}))_{t}R_{M}=JR_{M}=JT_{N}$.
Hence $(J^{2}(J:J^{2}))_{t_{1}}=J$ and therefore $T$ is Clifford
$t$-regular. \qed

{\bf Proof of Main Theorem.} On account of Lemma~\ref{sec:1.0-2}, we need only prove the ``only if"
assertion.

\texttt{Uniqueness}: Suppose there exist $J, F\in\{T\}\cup
\Gamma(I)$ such that $\overline{J}=\overline{F}$. Then there is
$0\not=q\in\qf(R)$ such that
$qJ=F=(F^{2})_{t}=(q^{2}J^{2})_{t}=q^{2}(J^{2})_{t}=q^{2}J$. So
$J=qJ=F$.

\texttt{Existence}: Let $J:=(I(T:I))_{t_{1}}=(I(T:I))_{t}$, a $t$-ideal of $T$ and a fractional
$t$-ideal of $R$  (by Lemma~\ref{sec:1.0-1}). Since $(I^{2})_{t}=qI$ for some  $0\not=q\in
\qf(R)$, then $(I:I^{2})=(I:(I^{2})_{t})=(I:qI)=q^{-1}(I:I)=q^{-1}T$. Hence
$J=(I(T:I))_{t}=(I(I:I^{2}))_{t}=(q^{-1}I)_{t}=q^{-1}I$. Therefore $\overline{J}=\overline{I}$. Now
$R$ is Clifford $t$-regular, then $I=(I^{2}(I:I^{2}))_{t}=(IJ)_{t}$. So $(I:I)\subseteq
(J:J)\subseteq ((IJ)_{t}:(IJ)_{t})=(I:I)$, whence $(J:J)=T$. Moreover,
$(T:J)=((I:I):J)=(I:IJ)=(I:(IJ)_{t})=(I:I)=T$. It follows that $(J:J^{2})=((J:J):J)=(T:J)=T$.
Consequently, $J=(J^{2}(J:J^{2}))_{t}=(J^{2})_{t}$ and thus $J$ is a fractional $t$-idempotent
$t$-ideal of $R$, and hence  a $t$-idempotent
$t$-ideal of $T$ by Lemma~\ref{sec:1.0-1}.\\
Now assume that $J\not=T$. Then we shall prove that $J\in\Gamma(I)$. By \cite[Theorem 3.2]{KM2} and Lemma~\ref{sec:1.7}, $T$ is a Krull-type domain. Then $J$ is contained in a finite number of $t$-maximal
ideals of $T$, say, $N_{1}, \dots, N_{r}$. By Lemma~\ref{sec:1.1}, for each $i\in \{1, \dots, r\}$,
$JT_{N_{i}}$ is an idempotent prime ideal of the valuation domain $T_{N_{i}}$. So
$JT_{N_{i}}=Q_{i}T_{N_{i}}$ for some prime ideal $Q_{i}\subseteq N_{i}$ of $T$; moreover,
$Q_{i}$ is minimal over $J$.  Then $J=\bigcap_{N\in\Max_{t}(T)}JT_{N}=(\bigcap_{1\leq
i\leq r}Q_{i}T_{N_{i}})\cap (\bigcap_{N'}JT_{N'})$, where $N'$ ranges over the  $t$-maximal
ideals of $T$ which do not contain $J$. The contraction to $T$ of both sides yields
$J=\bigcap_{1\leq i\leq r}Q_{i}$. One may assume the $Q_{i}$'s to be distinct. Since $JT_{N_{i}}$
is idempotent, $Q_{i}T_{N_{i}}=Q_{i}^{2}T_{N_{i}}$. We claim that $N_{i}$ is the unique $t$-maximal
ideal of $T$ containing $Q_{i}$. Otherwise, if $Q_{i}\subseteq N_{j}$ for some $j\not=i$, then
$Q_{i}$ and $Q_{j}$ are $t$-prime ideals \cite[Corollary 2.47]{Kg} contained in the same
$t$-maximal ideal $N_{j}$ in the PVMD $T$. So $Q_{i}$ and $Q_{j}$ are comparable under inclusion.
By minimality, we get $Q_{i}=Q_{j}$, absurd. It follows that
$Q_{i}=\bigcap_{j}Q_{i}T_{N_{j}}=Q_{i}T_{N_{i}}=Q_{i}^{2}T_{N_{i}}=\bigcap_{j}Q_{i}^{2}T_{N_{j}}=(Q_{i}^{2})_{t}$,
as desired. Finally, for each $i\in \{1, \dots r\}$, $T\subseteq (T:Q_{i})\subseteq (T:J)=(J:J)=T$.
Then $(T:Q_{i})=(Q_{i}:Q_{i})=T$. By \cite[Lemma 3.6]{KM2}, $Q_{i}$ is a $t$-maximal ideal of $T$,
completing the proof of the first statement.

Next, we describe the constituents groups $G_{\overline{J}}$ of $S_{t}(R)$. We write $G_{J}$
instead of $G_{\overline{J}}$, since we can always choose $J$ to be the unique fractional
$t$-idempotent $t$-ideal (of $R$) representing $\overline{J}$. We shall use $[L]$ to denote the
elements of the class group $\Cl(T)$. Notice that for any two common $t$-ideals $L, L'$ of $R$ and
$T$, we have $[L]=[L']$ if and only if $\overline{L}=\overline{L'}$ if and only if $L=xL'$ for some
$0\not=x\in\qf(R)=\qf(T)$.

\1 Assume that $J=T$. Let $[L]\in \Cl(T)$, where $L$ is a
$t$-invertible $t$-ideal of $T$. Then $T\subseteq (L:L)\subseteq
((L(T:L))_{t_{1}}:(L(T:L))_{t_{1}})=(T:T)=T$. Hence $(L:L)=T=(J:J)$.
We obtain, via Lemma~\ref{sec:1.0-1},
$(JL(L:L^{2}))_{t}=(JL(T:L))_{t}=(L(T:L))_{t}=(L(T:L))_{t_{1}}=T=J$.
Whence $\overline{L}\in G_{J}$. Conversely, let $\overline{L}\in
G_{J}$ for some $t$-ideal $L$ of $R$. By Lemma~\ref{sec:1.3-0},
$(L:L)=(J:J)$ and $(JL(L:L^{2}))_{t}=J$. By Lemma~\ref{sec:1.7},
$L$ is a $t$-ideal of $(L:L)=(J:J)=T$. Moreover, via
Lemma~\ref{sec:1.0-1},
$(L(T:L))_{t_{1}}=(L(T:L))_{t}=(L(L:L^{2}))_{t}=(JL(L:L^{2}))_{t}=J=T$.
Therefore $L$ is a $t$-invertible $t$-ideal of $T$ and thus $[L]\in
\Cl(T)$. Consequently, $G_{J}\cong\Cl(T)$.

\2 Assume $J=\bigcap _{1\leq i\leq r}Q_{i}$, where the $Q_{i}$'s are
distinct $t$-idempotent $t$-maximal ideals of $T$.

\begin{claim} $(T:J)=(J:J)=T$.\end{claim}

Indeed, for each $i$, we have $(Q_{i}:Q_{i})=T$ and hence
$(T:Q_{i})=((Q_{i}:Q_{i}):Q_{i})=(Q_{i}:Q_{i}^{2})=(Q_{i}:(Q_{i}^{2})_{t})=(Q_{i}:Q_{i})=T$.
Obviously, $\bigcap _{1\leq i\leq r}Q_{i}$ is an irredundant
intersection, so $(T:J)$ is a ring by \cite[Proposition
3.13]{HKLM2}. Further since $\{Q{i}\}_{1\leq i\leq r}$ equals the
set of minimal primes of $J$ in the PVMD $T$, then, by \cite[Theorem
4.5]{HKLM2}, $(T:J)=(\bigcap_{1\leq i\leq r}T_{Q_{i}})\cap
(\bigcap_{N'}T_{N'})$, where $N'$ ranges over the $t$-maximal ideals
of $T$ which do not contain $J$. Hence $(T:J)=T$, as claimed.

\begin{claim}$\phi$ is well-defined and injective.\end{claim}

Let $[L]\in \Cl(T)$ for some $t$-invertible $t$-ideal $L$ of $T$, that is,
$(L(T:L))_{t}=(L(T:L))_{t_{1}}=T$.  The homomorphism $\phi$ is given by
$\phi([L])=\overline{(LJ)_{t}}$. We have $(J:J)\subseteq (LJ:LJ)\subseteq ((LJ)_{t}:(LJ)_{t})$.
Conversely, let $x\in ((LJ)_{t}:(LJ)_{t})$. Then $x(LJ)_{t}\subseteq (LJ)_{t}$, hence
$x(LJ)_{t}(T:L)\subseteq (LJ)_{t}(T:L)$. So
$xJ=xJT=x(JT)_{t}=x(J(L(T:L))_{t})_{t}=x(JL(T:L))_{t}\subseteq
(JL(T:L))_{t}=(J(L(T:L))_{t})_{t}=J$. Therefore $x\in (J:J)$ and hence
$T=(J:J)=((LJ)_{t}:(LJ)_{t})$. Moreover,
$((LJ)_{t}(T:(LJ)_{t}))_{t}=(JL(T:JL))_{t}=(JL((T:J):L))_{t}=(JL(T:L))_{t}=(J(L(T:L))_{t})_{t}=JT=J$.
By Lemma~\ref{sec:1.3-0}, $\overline{(LJ)_{t}}\in G_{J}$ and thus $\phi$ is well-defined.

Now, let $[A]$ and $[B]$ in $\Cl(T)$ with $\overline{(AJ)_{t}}= \overline{(BJ)_{t}}$. So there
exists  $x\not=0\in \qf(R)=\qf(T)$ such that $(AJ)_{t}=x(BJ)_{t}$. Since $A$ and $B$ are
$t$-invertible $t$-ideals of $T$, then $A$ and $B$ are $v$-ideals of $T$. Further
$(T:A)=((T:J):A)=(T:JA)=(T:(AJ)_{t})=(T:x(BJ)_{t})=x^{-1}(T:(BJ)_{t})=x^{-1}(T:JB)=
x^{-1}((T:J):B)=x^{-1}(T:B)$. Hence $A=A_{v_{1}}=xB_{v_{1}}=xB$, whence $[A]=[B]$, proving that
$\phi$ is injective.

\begin{claim}
Let $Q$ be a $t$-idempotent $t$-maximal ideal of $T$ and $L$ a $t$-ideal of $T$ such that
${\overline{LT_{Q}}}\in G_{QT_{Q}}$. Then there exists a $t$-ideal $A$ of $T$ such that
${\overline{LT_{Q}}}={\overline{ AT_{Q}}}$, $Q$ is the unique $t$-maximal ideal of $T$ containing
$A$, and $\overline{A}\in G_{Q}$ in $S_{t}(T)$.
\end{claim}

We may assume $LT_{Q}\not=T_{Q}$, i.e., $L\subseteq Q$.
By \cite[Theorem 3.2]{KM2} and Lemma~\ref{sec:1.7}, $T$ is a Krull-type
domain. So let $\{Q, Q_{1}, ..., Q_{s}\}$ be the set of all
$t$-maximal ideals of $T$ containing $L$. Since $\{Q\wedge
Q_{i}\}_{1\leq i\leq s}$ is linearly ordered, we may assume that
$Q\wedge Q_{i}\subseteq P:= Q\wedge Q_{1}$ for each $i$.
Necessarily, $PT_{Q}\subsetneqq QT_{Q}$. On the other hand,
Lemma~\ref{sec:1.3-0} yields $(LT_{Q}:LT_{Q})=(QT_{Q}:QT_{Q})$ and
$LT_{Q}(QT_{Q}:LT_{Q})=QT_{Q}$; notice that in the valuation domain
$T_{Q}$, the $t$- and trivial operations coincide. Therefore there
exists $x\in (QT_{Q}:LT_{Q})$ such that $xLT_{Q}\nsubseteqq PT_{Q}$.
As $T_{Q}$ is valuation, $PT_{Q}\subsetneqq xLT_{Q}\subseteq
QT_{Q}$. Let $A:= xLT_{Q}\cap T$. Clearly, $P\subseteq A\subseteq Q$
and ${\overline{AT_{Q}}}={\overline{LT_{Q}}}$ (since
$AT_{Q}=xLT_{Q}$). Now assume there is $N\in \Max_{t}(T)$ such that
$A\subseteq N$ and $N\not= Q$. Since $L\subseteq Q_{1}$, then
$L\subseteq P\subseteq A\subseteq N$. Hence $N=Q_{i}$ for some $i\in
\{1, ..., s\}$, whence $A\subseteq Q\wedge Q_{i}\subseteq P$. So
$A=P$ and $PT_{Q}=xLT_{Q}$, absurd. Consequently, $Q$ is
the unique $t$-maximal ideal of $T$ containing $A$. Finally, one can
assume $A$ to be a $t$-ideal since $A\subseteq A_{t_{1}}\subseteq Q$
and $A_{t_{1}}T_{Q}=AT_{Q}$ by \cite[Lemma 3.3]{KM2}.

Next we show that $\overline{A}\in G_{Q}$ via Lemma~\ref{sec:1.3-0}. Since
$\overline{AT_{Q}}\in G_{QT_{Q}}$, then
$(AT_{Q}:AT_{Q})=(QT_{Q}:QT_{Q})=T_{Q}$ and
$AT_{Q}(T_{Q}:AT_{Q})=QT_{Q}$. Now, $(A:A)\subseteq T$ since
$A$ is an integral ideal of $T$. Conversely, we readily have $(A:A)\subseteq
(AT_{Q}:AT_{Q})=T_{Q}$ and $(A:A)\subseteq T_{N}$ for each $t$-maximal ideal $N\not=Q$ of $T$.
So that $(A:A)\subseteq T$ and hence $(A:A)=T$. Let $x\in
(T:A)$. Then $x\in (T_{Q}:AT_{Q})=(QT_{Q}:AT_{Q})$ since $AT_{Q}(T_{Q}:AT_{Q})=QT_{Q}$.
Therefore $xA\subseteq xAT_{Q}\subseteq QT_{Q}$ and hence
$xA\subseteq QT_{Q}\cap T=Q$, i.e., $x\in (Q:A)$. It follows that $(T:A)=(Q:A)$ and thus
 $A(T:A)=A(Q:A)$. Consequently, $(A(T:A))_{t_{1}}\subseteq Q$. Now, by the first statement of the theorem applied to $T$,
there exists a unique $t$-idempotent $t$-ideal $E$ of $T$ such that $\overline{A}\in G_{E}$ with
either $E=S$ for some fractional $t$-linked overring $S$ of $T$ or $E=\bigcap _{1\leq i\leq
s}N_{i}$, where the $N_{i}$'s are distinct $t$-idempotent $t$-maximal ideals of $S$. If $E=S$, then
$\overline{A}\in G_{S}$ implies that $(A:A)=(S:S)=S$ and $(A(S:A))_{t_{1}}=S$. So
$T=(A:A)=S=(A(S:A))_{t_{1}}\subseteq(A(T:A))_{t_{1}}\subseteq Q$, absurd. Hence, necessarily,
$E=\bigcap _{1\leq i\leq s}N_{i}$. It follows that $T=(A:A)=(E:E)=S$ (by the first claim) and
$(A(T:A))_{t_{1}}=E$. Therefore $A\subseteq E\subseteq N_{i}$ for each $i$, hence $E=Q$, the unique
$t$-maximal ideal of $T$ containing $A$. Thus, $\overline{A}\in G_{Q}$, proving the claim.

\begin{claim}$\psi$ is well-defined and surjective.\end{claim}

Let $\overline{L}\in G_{J}$ for some $t$-ideal $L$ of $R$. Notice
that $L$ is also a $t$-ideal of $(L:L)=(J:J)=T$. The homomorphism
$\psi$ is given by
$\psi(\overline{L})=(\overline{LT_{Q_{i}}})_{1\leq i\leq r}$. We
prove that $\psi$ is well-defined via a combination of \cite[Lemma
3.3]{KM2}, Lemma~\ref{sec:1.0-1}, and Lemma~\ref{sec:1.3-0}. Indeed
it suffices to show that
$(LT_{Q_{i}}:LT_{Q_{i}})=(Q_{i}T_{Q_{i}}:Q_{i}T_{Q_{i}})=T_{Q_{i}}$
and $LT_{Q_{i}}(T_{Q_{i}}:LT_{Q_{i}})=Q_{i}T_{Q_{i}}$. Let $i\in
\{1, ..., r\}$ and $x\in (LT_{Q_{i}}:LT_{Q_{i}})$. Then
$xLT_{Q_{i}}\subseteq LT_{Q_{i}}$ and hence
$xL(T:L)T_{Q_{i}}\subseteq L(T:L)T_{Q_{i}}$. So that $JT_{Q_{i}}=
(L(T:L))_{t}T_{Q_{i}}=(L(T:L))_{t_{1}}T_{Q_{i}}=(L(T:L))T_{Q_{i}}$.
Further $JT_{Q_{i}}=Q_{i}T_{Q_{i}}$. It follows that
$xQ_{i}T_{Q_{i}}\subseteq Q_{i}T_{Q_{i}}$, as desired. On the other
hand, we have $Q_{i}T_{Q_{i}}=JT_{Q_{i}}=(L(T:L))T_{Q_{i}}\subseteq
LT_{Q_{i}}(T_{Q_{i}}:LT_{Q_{i}})\subseteq T_{Q_{i}}$. The last
containment is necessarily strict. Otherwise $LT_{Q_{i}}=aT_{Q_{i}}$
for some $0\not=a\in L$. Therefore $L=(LJ)_{t}$ implies
$aT_{Q_{i}}=LT_{Q_{i}}=(LJ)_{t}T_{Q_{i}}=(LJ)_{t_{1}}T_{Q_{i}}=LJT_{Q_{i}}=aQ_{i}T_{Q_{i}}$,
absurd. Consequently
$LT_{Q_{i}}(T_{Q_{i}}:LT_{Q_{i}})=Q_{i}T_{Q_{i}}$. So $\psi$ is
well-defined.

Next we show that $\psi$ is surjective. Let
$(\overline{L_{i}T_{Q_{i}}})_{1\leq i\leq r}\in \prod
G_{Q_{i}T_{Q_{i}}}$. By the above claim, for each $i$, there
exists a $t$-ideal $A_{i}$ of $T$ such that
${\overline{L_{i}T_{Q_{i}}}}={\overline{ A_{i}T_{Q_{i}}}}$, $\overline{A_{i}}\in G_{Q_{i}}$, and
$Q_{i}$ is the unique $t$-maximal ideal of $T$ containing $A_{i}$. Set $A:=(A_{1}A_{2}\dots A_{r}J)_{t_{1}}=(A_{1}A_{2}\dots
A_{r}J)_{t}$. Let $j\in \{1, \dots, r\}$. By \cite[Lemma 3.3]{KM2},
$AT_{Q_{j}}=(A_{1}A_{2}\dots A_{r}J)T_{Q_{j}}=A_{j}Q_{j}T_{Q_{j}}$
since $JT_{Q_{j}}=Q_{j}T_{Q_{j}}$ and $A_{i}T_{Q_{j}}=T_{Q_{j}}$ for
each $i\not =j$. So
$\overline{AT_{Q_{j}}}=\overline{A_{j}Q_{j}T_{Q_{j}}}={\overline{A_{j}T_{Q_{j}}}}\ {\overline{Q_{j}T_{Q_{j}}}}=\overline{A_{j}T_{Q_{j}}}=\overline{L_{j}T_{Q_{j}}}$. Therefore
$\psi(\overline{A})=(\overline{L_{i}T_{Q_{i}}})_{1\leq i\leq r}$.

Next we show that $\overline{A}\in G_{J}$. First notice that $Q_{1}, \dots, Q_{r}$ are the only
$t$-maximal ideals of $T$ containing $A$. For, let $Q$ be a $t$-maximal ideal of $T$ such that
$A\subseteq Q$. Then either $J\subseteq Q$ or $A_{i}\subseteq Q$ for some $i$. In both cases,
$Q=Q_{j}$ for some $j$, as desired. Now $A$ is an ideal of $T$, then $(J:J)=T\subseteq (A:A)$.
Conversely, for each $j$, $A_{j}Q_{j}T_{Q_{j}}=a_{j}A_{j}T_{Q_{j}}$, for some nonzero
$a_{j}\in\qf(T)$, since $\overline{A_{j}T_{Q_{j}}}\in G_{Q_{j}T_{Q_{j}}}$. So $(A:A)\subseteq
(AT_{Q_{j}}:AT_{Q_{j}})=(A_{j}Q_{j}T_{Q_{j}}:A_{j}Q_{j}T_{Q_{j}})
=(a_{j}A_{j}T_{Q_{j}}:a_{j}A_{j}T_{Q_{j}})=(A_{j}T_{Q_{j}}:A_{j}T_{Q_{j}})=T_{Q_{j}}$. Further, for
each $N\in\overline{\Max_{t}}(T, A)$, we clearly have $(A:A)\subseteq T_{N}$. It follows that
$(A:A)\subseteq T=(J:J)$ and hence $(A:A)=T=(J:J)$. Next we prove that $(A(T:A))_{t}=J$. We have
$(A_{i}:A_{i})=T$ and $(A_{i}(T:A_{i}))_{t}=(A_{i}(T:A_{i}))_{t_{1}}=Q_{i}$, for each $i$, since
$\overline{A_{i}}\in G_{Q_{i}}$. Let $j\in\{1,\dots, r\}$ and set $F_{j}:=\prod_{i\not
=j}A_{i}$. Clearly $A=(JF_{j}A_{j})_{t}$. Now, since $A\subseteq A_{j}$, then
$(A(T:A_{j}))_{t}\subseteq (A(T:A))_{t}$. However
$(A(T:A_{j}))_{t}=((JF_{j}A_{j})_{t}(T:A_{j}))_{t}=(JF_{j}A_{j}(T:A_{j}))_{t}=(JF_{j}(A_{j}(T:A_{j}))_{t})_{t}=(JF_{j}Q_{j})_{t}$.
Hence $(JF_{j}Q_{j})_{t}\subseteq (A(T:A))_{t}$. So
$JQ_{j}T_{Q_{j}}=JF_{j}Q_{j}T_{Q_{j}}=(JF_{j}Q_{j})_{t}T_{Q_{j}}\subseteq (A(T:A))_{t}T_{Q_{j}}$
since $F_{j}T_{Q_{j}}=T_{Q_{j}}$. Then
$JT_{Q_{j}}=Q_{j}T_{Q_{j}}=Q_{j}^{2}T_{Q_{j}}=JQ_{j}T_{Q_{j}}\subseteq (A(T:A))_{t}T_{Q_{j}}$.
Since $\Max_{t}(T, A)=\Max_{t}(T, J)$ and $\overline{\Max_{t}}(T, A)=\overline{\Max_{t}}(T,
J)$, it follows that $J\subseteq (A(T:A))_{t}$. Conversely, let $j\in \{1, ..., r\}$. By the proof
of the third claim, $(T:A_{j})=(Q_{j}:A_{j})$. Then
$(T:A)=(T:(JF_{j}A_{j})_{t})=(T:JF_{j}A_{j})=((T:A_{j}):JF_{j})=((Q_{j}:A_{j}):JF_{j})=(Q_{j}:JF_{j}A_{j})=(Q_{j}:(JF_{j}A_{j})_{t})=(Q_{j}:A)$.
Hence $(T:A)=\bigcap_{1\leq j\leq r}(Q_{j}:A)=((\bigcap_{1\leq j\leq r}Q_{j}):A)=(J:A)\subseteq
(T:A)$. So $(T:A)=(J:A)$. Therefore $(A(T:A))_{t}=(A(J:A))_{t}\subseteq J$. Consequently,
$J=(A(T:A))_{t}$ and thus $\overline{A}\in G_{J}$, as desired.

\begin{claim} $\Ima(\phi)=\Ker(\psi)$. \end{claim}
Indeed, let $[A]\in \Cl(T)$ for some $t$-invertible $t$-ideal $A$ of
$T$. Then there exists a finitely generated ideal $B$ of $T$ such
that $A=B_{v_{1}}=B_{t_{1}}$. Hence
$\psi(\phi([A]))=\psi(\overline{(AJ)_{t}})=(\overline{(AJ)_{t}T_{Q_{i}}})_{1\leq
i\leq r}$. For each $i$, we have
$(AJ)_{t}T_{Q_{i}}=(AJ)_{t_{1}}T_{Q_{i}}=(BJ)_{t_{1}}T_{Q_{i}}=BJT_{Q_{i}}=BQ_{i}T_{Q_{i}}=b_{i}Q_{i}T_{Q_{i}}$
for some nonzero $b_{i}\in B$. Then
$\overline{(AJ)_{t}T_{Q_{i}}}=\overline{Q_{i}T_{Q_{i}}}$ in
$G_{Q_{i}T_{Q_{i}}}$. It follows that $\Ima(\phi)\subseteq
\Ker(\psi)$.

Conversely, let $\overline{L}\in G_{J}$ such that
$\overline{LT_{Q_{i}}}=\overline{Q_{i}T_{Q_{i}}}$ for each $i\in
\{1, \dots r\}$, that is, there exists $a_{i}\not=~0\in \qf(T)$ such
that $a_{i}Q_{i}T_{Q_{i}}=LT_{Q_{i}}\subseteq T_{Q_{i}}$. Then
$a_{i}\in
(T_{Q_{i}}:Q_{i}T_{Q_{i}})=(Q_{i}T_{Q_{i}}:Q_{i}^{2}T_{Q_{i}})=(Q_{i}T_{Q_{i}}:Q_{i}T_{Q_{i}})=T_{Q_{i}}$
for each $i$. Let $B:=\sum_{1\leq k\leq r} Ta_{i}$ and
$A:=B_{t_{1}}$. Clearly, $A$ is a fractional $t$-invertible
$t$-ideal of $T$, i.e., $[A]\in \Cl(T)$. Further, for each $i$,
$(AJ)_{t}T_{Q_{i}}=(AJ)_{t_{1}}T_{Q_{i}}=(BJ)_{t_{1}}T_{Q_{i}}=BJT_{Q_{i}}=BQ_{i}T_{Q_{i}}=a_{k}Q_{i}T_{Q_{i}}$
for some  $a_{k}\ (\not= 0)$, hence
$\overline{(AJ)_{t}T_{Q_{i}}}=\overline{LT_{Q_{i}}}$. Therefore
$\phi([A])=\overline{L}$. Hence $\Ker(\psi)\subseteq \Ima(\phi)$, as
desired.

Consequently, the sequence is exact, completing the proof of the theorem. \qed

A domain $R$ is said to be strongly $t$-discrete if it has no
$t$-idempotent $t$-prime ideals, i.e., for every  $t$-prime ideal
  $P$ of $R$, $(P^{2})_{t}\subsetneq P$.

\begin{corollary}\label{sec:1.8} Let $R$ be a Krull-type domain
which is strongly $t$-discrete. Then $S_{t}(R)\cong\bigvee_{T} \Cl(T)$, where $T$ ranges over the
set of fractional $t$-linked overrings of $R$.
\end{corollary}

\proof  Recall first the fact that every fractional $t$-linked overring $T$ of
$R$ has the form $T=(I:I)$ for some $t$-ideal $I$ of $R$ such that $\overline{I}$ is an idempotent
of $S_{t}(R)$; precisely, $I:=aT$, for some $0\not=a\in (R:T)$, with
$(I^{2})_{t}=(I^{2})_{t_{1}}=a^{2}T=aI$. Now, Lemma~\ref{sec:1.10} forces each $T$ to be strongly $t$-discrete and then Theorem~\ref{sec:1.2} leads to the conclusion. \qed

Since in a Pr\"ufer domain the $t$-operation coincides with the
trivial operation, we recover Bazzoni's results on Pr\"ufer domains
of finite character.

\begin{corollary}[{\cite[Theorem 3.1]{Ba2} \& \cite[Theorem 3.5]{Ba3}}] \label{sec:1.9}
Let $R$ be a Pr\"ufer domain of finite character. Then $\overline{J}$ is an idempotent of $S(R)$ if
and only if there exists a unique nonzero idempotent fractional ideal $L$ such that $\overline{J}=
\overline{L}$ and $L$ satisfies one of the
following two conditions:\\
(1) $L:=D$ where $D$ is a fractional overring of $R$, or\\
(2) $L:=P_{1}.P_{2}...P_{n}D$, where each $P_{i}$ is a nonzero
idempotent prime ideal of $R$, and $D$ is a fractional overring of
$R$. Moreover, the following sequence is exact
$$0\longrightarrow \Cl(D)\longrightarrow G_{L}\longrightarrow
\prod_{1\leq i\leq r}G_{P_{i}R_{P_{i}}}\longrightarrow 0.$$
\end{corollary}

\proof     The result follows readily from Theorem~\ref{sec:1.2} since $T$ is
flat over $R$ and every prime ideal $Q$ of $T$ is of the form $Q=PT$ for some prime $P$ of $R$, and
$Q$ is idempotent if and only if so is $P$. \qed

\section{Examples}\label{sec:2}

One can develop numerous illustrative examples via Theorem~\ref{sec:1.2} and Corollary~\ref{sec:1.8}. We'll provide two families of such examples by means of polynomial rings over valuation domains. For this purpose, we first state the following lemma.

\begin{lemma}\label{sec:2.1}
Let $V$ be a nontrivial valuation domain, $X$ an indeterminate
over $V$, and $R:=V[X]$. Then the following statements hold:\\
\1 $R$ is a Krull-type domain which is not Pr\"ufer.\\
\2 Every fractional $t$-linked overring of $R$ has the form $V_{p}[X]$ for some nonzero $p\in\Spec(V)$.\\
\3 Every $t$-idempotent $t$-prime ideal of $R$ has the form $p[X]$
for some idempotent $p\in\Spec(V)$.
\end{lemma}

\proof     \1  The notion of PVMD is stable
under adjunction of indeterminates \cite{HMM}. So $R$ is a PVMD and
has finite $t$-character by \cite[Proposition 4.2]{KM2}, as desired.
Further, $R$ is not Pr\"ufer since nor is any polynomial ring over a
nontrivial integral domain.\\
\2 Let $p\not=0\in\Spec(V)$, then $V_{p}[X]$ is a fractional
$t$-linked overring of $R$. Indeed, let $S:=V\setminus p$ and let $J$ be a finitely generated ideal of $R$ such that $J^{-1}=R$. We have
$(V_{p}[X]:JV_{p}[X])=(S^{-1}(V[X]):S^{-1}J)=S^{-1}(R:J)=S^{-1}R=V_{p}[X]$.
Hence $V_{p}[X]$ is $t$-linked over $R$. Now suppose $p$ is not maximal. Since $V$ is a conducive domain (since valuation), then $(V:V_{p})=p$ \cite{DF}. Hence $(V[X]:V_{p}[X])=(V:V_{p})[X]=p[X]$. It follows that $V_{p}[X]$ is a
fractional overring of $R$. If $p$ is maximal, then $V_{p}[X]=V[X]=R$ is trivially a fractional overring of $R$. Next let $T$ be a
fractional $t$-linked overring of $R$, $0\not=a\in (R:T)$, and
$I:=aT$. By Lemma~\ref{sec:1.0-1}, $I$ is a common $t$-ideal of both
$R$ and $T$. Set $A:= I\cap V$. If $A\not= 0$, then $A$ is a
$t$-ideal of $V$ and hence $I = A[X]$; if $A=(0)$, then $I=fB[X]$
where $f\not=0\in \qf(V)[X]$ and $B$ is a $t$-ideal of $V$ \cite{Q}.
If $I=A[X]$, then $T=(I:I)=(A[X]:A[X])=(A:A)[X]$; and if $I=fB[X]$,
then $T=(I:I)=(fB[X]:fB[X])=(B[X]:B[X])=(B:B)[X]$. Moreover, $(A:A)$
and $(B:B)$ are overrings (and hence localizations) of $V$.
Therefore, in both cases, $T=V_{p}[X]$ for some nonzero prime ideal
$p$ of $V$, as desired.\\
\3 Let $p$ be an idempotent prime ideal of $V$. Then
$((p[X])^{2})_{t}=(p^{2}[X])_{t}=(p[X])_{t}=p[X]$, recall here that the $t$-operation with respect
to $V$ is trivial. Next let $P$  be a $t$-idempotent $t$-prime ideal of $R$ and $p:=P\cap V$.
Assume that $p=(0)$ and set $S:=V\setminus\{0\}$. Then, by \cite[Lemma 2.6]{KM2}, $S^{-1}P$ is an
idempotent (nonzero) ideal of the PID $S^{-1}R=\qf(V)[X]$, absurd. It follows that $p\not=(0)$.
Since $V$ is integrally closed, then $P=p[X]$ \cite{Q}. Moreover $p[X]=(p^{2}[X])_{t}=p^{2}[X]$,
hence $p=p^{2}$, as desired. \qed

\begin{example}\label{sec:2.2}\rm
Let $n$ be an integer $\geq 1$. Let $V$ be an $n$-dimensional strongly discrete valuation domain
and let $(0)\subset p_{1}\subset p_{2}\subset ...\subset p_{n}$ denote the chain of its prime
ideals. Let $R:=V[X]$, a Krull-type domain. A combination of Lemma~\ref{sec:2.1} and
Corollary~\ref{sec:1.8} yields $S_{t}(R)\cong\bigvee_{1\leq i\leq n} \Cl(V_{p_{i}}[X])$. Moreover
$\Cl(V_{p_{i}}[X])=\Cl(V_{p_{i}})=0$, so that ${\mathcal S_t}(R)$ is a disjoint union of $n$ groups
all of them are trivial. Precisely, ${\mathcal S_t}(R)=\{V_{p_{1}}[X], V_{p_{2}}[X], ...,
V_{p_{n}}[X]\}$ where, for each $i$, $\overline{V_{p_{i}}[X]}$ is identified with $V_{p_{i}}[X]$
(due to the uniqueness stated by Theorem~\ref{sec:1.2}). \qed\end{example}

\begin{example}\label{sec:2.3}\rm
Let $V$ be a one-dimensional valuation domain with idempotent maximal ideal $M$ and $X$ an
indeterminate over $V$. Let $R:=V[X]$, a Krull-type domain. By
Theorem~\ref{sec:1.2} and Lemma~\ref{sec:2.1}, we obtain $S_{t}(R)=\Cl(R)\vee G_{M[X]}$.
$\Cl(R)=\Cl(V[X])=\Cl(V)=0$. So $S_{t}(R)=\{\overline{R}\}\vee \{\overline{I}\ |\ I\
t\mbox{-ideal of}\ R\ \mbox{with}\ (II^{-1})_{t}=M[X]\}.$ \qed\end{example}


\end{document}